\documentclass{amsart}
\usepackage{amsmath,amsfonts}
\usepackage{latexsym}
\title[C$^*$-algebra-valued-symbol pseudos]{C$^*$-algebra-valued-symbol 
pseudodifferential operators: 
abstract characterizations}
\author[S. T. Melo \and M. I. Merklen]{Severino T. Melo 
\and Marcela I. Merklen}
\date{}
\newtheorem{thm}{Theorem}

\newtheorem{lem}{Lemma}
\newtheorem{cor}{Corollary}
\newtheorem{df}{Definition}
\begin{document}
\newcommand{\re}{{\mathbb R}}
\newcommand{\C}{{\mathbb C}}
\newcommand{\N}{{\mathbb N}}
\newcommand{\rn}{{\mathbb R}^{n}}
\newcommand{\rtn}{{\mathbb R}^{2n}}
\newcommand{\op}{operator}
\newcommand{\ops}{operators}
\newcommand{\psd}{pseudo\-dif\-fer\-en\-tial}
\newcommand{\Psd}{Pseudo\-dif\-fer\-en\-tial}
\newcommand{\bC}{{\mathcal B}^C(\rtn)}
\newcommand{\bcn}{{\mathcal B}^C(\rn)}
\newcommand{\ben}{{\mathcal B}^*(E_n)}
\newcommand{\sC}{{\mathcal S}^C(\rn)}
\newcommand{\bCd}{{\mathcal B}^C(\re^2)}
\newcommand{\sCt}{{\mathcal S}^C(\rtn)}
\newcommand{\cqd}{\hfill$\Box$}
\newcommand{\talg}{\stackrel{\mbox{{\tiny{\tt alg}}}}{\otimes}}
\newcommand{\du}{d\!u}
\newcommand{\dv}{d\!v}
\newcommand{\dx}{d\!x}
\newcommand{\dz}{d\!z}
\newcommand{\dy}{d\!y}
\newcommand{\deta}{d\!\eta}
\newcommand{\dxi}{d\!\xi}
\newcommand{\dzeta}{d\!\zeta}
\newcommand{\h}{{\mathcal H}}
\newcommand{\pf}{{\em Proof}: }

\begin{abstract}

Given a separable unital C$^*$-algebra $C$ with norm $||\cdot||$, let $E_n$ 
denote the Banach-space completion of the $C$-valued Schwartz space on $\rn$ 
with norm $||f||_2=||\langle f,f\rangle||^{1/2}$, 
$\langle f,g\rangle=\int f(x)^*g(x)\dx$. The assignment of the 
pseudodifferential operator $A=a(x,D)$ with $C$-valued symbol 
$a(x,\xi)$ to each smooth function with bounded derivatives $a\in\bC$ 
defines an injective mapping $O$, from $\bC$ to the set $\h$ of all 
operators with smooth orbit under the canonical action of the Heisenberg group 
on the algebra of all adjointable operators on the Hilbert module $E_n$. In 
this paper, we construct a left-inverse $S$ for $O$ and prove that $S$ is 
injective if $C$ is commutative. This generalizes Cordes' description \cite{C} 
of $\h$ in the scalar case. Combined with previous results of the second-named 
author, our main theorem implies that, given a skew-symmetric 
$n\times n$ matrix $J$, and if $C$ is commutative, then any $A\in\h$ which 
commutes with every  \psd\ \op\ with symbol $F(x+J\xi)$, 
$F\in {\mathcal B}^C(\rn)$, is a \psd\ \op\ with symbol $G(x-J\xi)$, for some 
$G\in {\mathcal B}^C(\rn)$. That was conjectured by Rieffel.
\end{abstract}
\maketitle

\begin{center}
\end{center}

\section{Introduction}

Let $C$ be a separable unital C$^*$-algebra with norm $||\cdot||$, and
let $\sC$ denote the set of all $C$-valued smooth functions on $\rn$ which, 
together with all their derivatives, are bounded by arbitrary negative 
powers of $|x|$, $x\in\rn$. We equip it with the $C$-valued inner-product 
\[
\langle f,g\rangle\,=\,\int f(x)^*g(x)\dx,
\]
which induces the norm $||f||_2=||\langle f,f\rangle||^{1/2}$, and denote by 
$E_n$ its Banach-space completion with this norm. The inner product 
$\langle\cdot,\cdot\rangle$ 
turns $E_n$ into a Hilbert module \cite{L}. The set of all (bounded) adjointable \ops\ 
on $E_n$ is denoted ${\mathcal B}^*(E_n)$. 

Let $\bC$ denote the set of all 
smooth bounded functions from $\rtn$ to $C$ whose derivatives of arbitrary 
order are also bounded. For each $a$ in $\bC$, a linear mapping from $\sC$ to 
itself is defined by the formula
\begin{equation}
\label{psddef}
(Au)(x)\,=\,\frac{1}{(2\pi)^{n/2}}\int e^{ix\cdot\xi}a(x,\xi)\hat u(\xi)d\xi,
\end{equation}
where $\hat u$ denotes the Fourier transform, 
\[
\hat u(\xi)=(2\pi)^{-n/2}\int e^{-iy\cdot\xi}u(y)\dy.
\] 
As usual, we denote $A=a(x,D)$. This operator extends to an element of 
${\mathcal B}^*(E_n)$ whose norm satisfies the following estimate. There 
exists a constant $k>0$ depending only on $n$ such that 
\begin{equation}
\label{cvm}
||A||\,\leq\,k \sup\{\,
||\partial_x^\alpha\partial_\xi^\beta a(x,\xi)||;
(x,\xi)\in\rtn\ \mbox{and}\ \alpha,\beta
\leq(1,\cdots,1)\,\}.
\end{equation}
This generalization of the Calder\'on-Vaillancourt Theorem 
\cite{CV} was proven by Mer\-klen \cite{M,Mt}, following ideas of Hwang
\cite{H}\ and Seiler \cite{S}. The case of 
$a(x,\xi)=F(x+J\xi)$, where $F\in{\mathcal B}^C(\rn)$ and $J$ is an 
$n\times n$ skew-symmetric matrix, had been proven earlier by Rieffel 
\cite[Corollary 4.7]{R}. 

The estimate \eqref{cvm} implies that the mapping 
\begin{equation}
\label{azz}
\rtn\ni(z,\zeta)\longmapsto  A_{z,\zeta}= T_{-z}M_{-\zeta}AM_{\zeta}
T_z\in{\mathcal B}^*(E_n)
\end{equation}
is smooth (i.e., $C^\infty$ with respect to the norm topology), where $T_z$
and $M_\zeta$ are defined by $T_zu(x)=u(x-z)$ and $M_{\zeta}u(x)=e^{i\zeta \cdot x}u(x)$, $u\in\sC$. 
That follows just like in the scalar case \cite[Chapter 8]{book}. 

\begin{df}
We call 
{\em Heisenberg smooth} an operator $A\in{\mathcal B}^*(E_n)$ for which the 
mapping \eqref{azz} is smooth, and denote by $\h$ the set of all such \ops. 
\end{df}

The elements of $\h$ are the smooth vectors for the action of the Heisenberg 
group on 
${\mathcal B}^*(E_n)$ given by the same formula as the standard one in the 
scalar case (i.e., when $C$ is the algebra $\C$ of complex numbers and then 
$E_n=L^2(\rn)$, and we denote $\sC$ and $\bC$ by ${\mathcal S}(\rn)$
and ${\mathcal B}(\rtn)$, respectively).

We therefore have a mapping
\begin{equation}
\label{o}
\begin{array}{rcl} 
O:\bC&\longrightarrow&\h\\
a\ &\longmapsto&O(a)=a(x,D).
\end{array}
\end{equation}
In the scalar case, it is well-known (this can be proven by a Schwartz-kernel 
argument) that if a \psd\ \op\ as in \eqref{psddef} vanishes on 
${\mathcal S}(\rn)$, then $a$ must be zero.  Let us show that this implies 
that $O$ is injective for arbitrary $C$. 

Given any complex-valued function $u$ defined 
on $\rn$, we denote by $\tilde u:\rn\to C$ the function defined by 
\begin{equation}
\label{tilde}
\tilde u(x)=u(x){\mathbf 1}_C,
\end{equation}
where ${\mathbf 1}_C$ denotes the identity of $C$. If $O(a)=0$, the fact that
$O(a)\tilde u=0$ for every $u\in{\mathcal S}(\rn)$ and the injectivity of $O$ 
in the scalar case imply that $(x,\xi)\mapsto \rho(a(x,\xi))$ vanishes 
identically, for every $\rho\in C^*$, the (Banach-space) dual of $C$. We then
get $a\equiv 0$, as we wanted. 

Our results in this paper can now be summarized in the following theorem, 
proven in Sections~\ref{s2} and~\ref{s3}.

\begin{thm}\label{t}
Let $C$ be a unital separable C$^*$-algebra.
There exists a linear mapping $S:\h\to\bC$ such that $S\circ O$ is the
identity operator. If $C$ is commutative, then $S$ is injective.
\end{thm}

Since an injective left-inverse is an inverse, we get:

\begin{cor}\label{c}If $C$ is commutative and an operator 
$A\in{\mathcal B}^*(E_n)$ is given, then the mapping defined in 
{\em \eqref{azz}} is smooth if and only if $A=a(x,D)$ for some 
$a\in\bC$.
\end{cor}

Theorem \ref{t} and Corollary \ref{c} were proven by Cordes \cite{C} in the 
scalar case. His construction \cite[Chapter 8]{book} of the left-inverse $S$  
works also in the general case, if only one is careful enough to avoid mentioning 
trace-class or Hilbert-Schmidt \ops. That is what we show in
Section~\ref{s2}. His proof that $S$ is injective, however, strongly depends on the 
fact that, when $C=\C$, $E_n=L^2(\rn)$ is a Hilbert space. In the general commutative
case, the lack of an orthonormal basis in $E_n$ can be bypassed by still reducing 
the problem to copies of $L^2(\rn)$, as shown at the beginning of Section~\ref{s3}. After this
reduction, we are then able to follow the steps of Cordes' proof. Crucial for this
strategy, our Lemma~\ref{lr} is essentially \cite[Lemma 2.4]{Ra}\ specialized to 
commutative C$^*$-algebras. In Section~\ref{s4}, we explain how Theorem~\ref{t} 
implies, in the commutative case, an abstract characterization, conjectured by 
Rieffel \cite{R}, of a certain class of C$^*$-algebra-valued-symbol \psd\ \ops.

The asumption of separability of $C$ is needed to justify several results
about vector-valued integration (see \cite[Ap\^endice]{Mt}, for example), which are used without
further comments throughout the text. 

\section{Left Inverse for $O$}\label{s2}

Given $f$ and $g$ functions from $\rn$ to $X$ ($X$ will be either
$C$ or $\C$), let $f\otimes g:\rtn\to X$ 
be defined by
\begin{equation}
\label{tensor}
f\otimes g(x,y)=f(x)g(y). 
\end{equation}
Given a vector space $V$ we denote by $V\talg V$ the algebraic tensor 
product of $V$ by itself. In case the elements of $V$ are functions 
from $\rn$ to $X$, $V\talg V$ is isomorphic to
the linear span of all function as in \eqref{tensor} with $f$ and $g$ in $V$.

\begin{lem} 
\label{otimes}
Given $A\in B^*(E_n)$ mapping $\sC$ to itself, there exists 
a unique operator $A\otimes I\in B^*(E_{2n})$ such that, for all $f$ an $g$ in 
$\sC$, 
\begin{equation}
\label{innt}
(A\otimes I)(f\otimes g)\,=\,Af\otimes g.
\end{equation}
\end{lem}

\pf
Let $L^2(\rn;C)$ denote the set of equivalence classes (for the
equality almost everywhere equivalence) of Borel measurable functions
$f:\rn\to C$ such that
\[
\int||f(x)||^2\dx\,<\,\infty.
\]
and let $||f||_{L^2}$ denote the square root of the integral above. 
$L^2(\rn;C)$ equipped with $||\cdot||_{L^2}$ is a Banach space, containing
$\sC$ as a dense subspace. It follows from the inequality
\[
||f||_2\,\leq\,||f||_{L^2},\ \mbox{for all}\  f\in\sC,
\]
that $L^2(\rn;C)$ embeds in $E_n$ as a $||\cdot||_2$-dense subspace. 

Let $\mbox{{\tt S}}_n$ denote the set of all simple measurable functions
from $\rn$ to $C$.  
It takes an elementary but messy argument to show that 
$\mbox{{\tt S}}_n\talg \mbox{{\tt S}}_n$ is $||\cdot||_{L^2}$-dense in 
$\mbox{{\tt S}}_{2n}$, which is dense in $L^2(\rtn;C)$.
Since $\mbox{{\tt S}}_n$ is dense in $L^2(\rn;C)$, it 
follows that $L^2(\rn;C)\talg L^2(\rn;C)$ is dense in $L^2(\rtn;C)$. 
Since $\sC$ is dense in  $L^2(\rn;C)$, it follows that $\sC\talg\sC$ is 
$||\cdot||_{L^2}$-dense in $L^2(\rtn;C)$, hence it is also 
$||\cdot||_2$-dense in $E_{2n}$. 

Let $\phi:C\to{\mathcal B}^*(E_n)$ be given by left multiplication
on $\sC$, and denote by $E_n\otimes_{\phi}E_n$ the interior tensor 
product (given by $\phi$) as defined in \cite[page 41]{L}. The fact that 
$\sC\talg\sC$ is dense in $E_{2n}$ allows us to identify 
$E_n\otimes_{\phi}E_n$ with $E_{2n}$ (notice that the space $N$ in 
\cite[Proposition 4.5]{L} consists only of $0$ in this case). 

Given $A\in{\mathcal B}^*(E_n)$, it now follows from the more general result 
around \cite[(4.6)]{L} that there exists a unique 
$A\otimes I\in{\mathcal B}^*(E_{2n})$ such that $A\otimes I(f\otimes g)=
Af\otimes g$ for all $f\otimes g\in E_n\talg E_n$. In particular, we get \eqref{innt} for all $f$ and 
$g$ in $\sC$. That \eqref{innt} uniquely determines $A\otimes I$ also 
follows from the fact that $\sC\talg\sC$ is dense in $E_{2n}$. \cqd

Let us denote by $\gamma_1(t)$ and $\gamma_2(t)$, respectively, the 
fundamental solutions of 
$(\partial_t+1)$ and $(\partial_t+1)^2$ given by:
\[
\gamma_1(t)=\left\{
\begin{array}{rl}
e^{-t},&\ \mbox{if}\ t\geq 0\\
0,&\ \mbox{if}\ t< 0
\end{array}
\right.\ \mbox{and}\
\gamma_2(t)=\left\{
\begin{array}{rl}
te^{-t},&\ \mbox{if}\ t\geq 0\\
0,&\ \mbox{if}\ t< 0
\end{array}
\right..
\]
We then define $u$ and $v$ in $L^2(\re)\cap L^1(\re)$ by
\begin{equation}
\label{defv}
v(\xi,\eta)=\gamma_1(\xi-\eta)/(1+i\xi)^2
\end{equation}
and
\begin{equation}
\label{defu}
u(x,\eta)=(1+\partial_{\eta})[(1-i\eta)^2\gamma_2(-x)
\gamma_2(-\eta)e^{ix\eta}].
\end{equation}
The following lemma can be proven exactly like in the scalar case
\cite[Section 8.3]{book}.

\begin{lem}\label{cordes} 
If $a$ and $b$ in $\bCd$ are such that
$(1+\partial_z)^2(1+\partial_\zeta)^2a(z,\zeta)=b(z,\zeta)$,
then we have, for all $(z,\zeta)\in\re^2$, 
\begin{equation}
\label{hoc}
a(z,\zeta)=\int_{\re^{3}}
\overline{u(x,\eta)}e^{ix\xi}b(x+z,\xi+\zeta)v(\xi,\eta)
\dxi\dx\deta.
\end{equation}
\end{lem}

We also omit the proof of the following lemma.

\begin{lem} 
\label{lemvl}
There exists a sequence $v_l$ in 
${\mathcal S}(\re)\talg{\mathcal S}(\re)$
such that $v_l\to v$ in $L^2(\re^2)$ and
\[
\lim_{l\to\infty}\,
\int_{\re^{3}}|u(x,\eta)|\cdot|v(\xi,\eta)-v_l(\xi,\eta)|\,\dxi\dx\deta
\,=\,0.
\]
\end{lem}

We are ready to define $S$ when $n=1$. Given $A\in\h$, let $B=f(0,0)$, where
$f:\re^2\to{\mathcal B}^*(E_1)$ denotes the smooth function 
\begin{equation}\label{bzz}
f(z,\zeta)\,=\,(1+\partial_z)^2(1+\partial_\zeta)^2A_{z,\zeta}\,, 
\end{equation}
with $A_{z,\zeta}$ as in \eqref{azz}. The group property allows one to show that
$f(z,\zeta)=B_{z,\zeta}$ for all $(z,\zeta)\in\re^2$. We then define 
\begin{equation}
\label{defs}
(SA)(z,\zeta)\ =\ \sqrt{2\pi}\,\langle\,\tilde u, (B_{z,\zeta}F^*\otimes I) 
\tilde v\,\rangle,
\end{equation}
where $F\in{\mathcal B}^*(E_1)$, $F^{*}=F^{-1}$, denotes the Fourier 
transform, and $\langle\cdot,\cdot\rangle$ denotes the inner product of
$E_2$. The meaning of $\tilde{\ }$ was defined in \eqref{tilde} and we are
regarding, as explained in the proof of Lemma~\ref{otimes}, $L^2(\re^2;C)$ as a
subspace of $E_2$.

It is not hard to see that $S$ maps $\h$ to $\bCd$ (this uses the inequality
$||A\otimes I||\leq||A||$, which follows from \cite[(4.6)]{L}). Given $a\in\bCd$, let $c=SOa$. 
To prove that $S\circ O$ is the identity on $\bCd$, it is enough
to show that
\[
\int_{\re^2}\,[a(z,\zeta)-c(z,\zeta)] f(z,\zeta)\,\dz\dzeta=0,
\ \ \mbox{for all}\ \ f\in{\mathcal S}(\re^2).
\]
Indeed, if this is the case, then 
$(z,\zeta)\mapsto\rho(a(z,\zeta)-c(z,\zeta))$ vanishes identically 
for all $\rho\in C^*$, and the equality $a=c$ will therefore hold. 

For each $l\in\N$, define 
$
c_l(z,\zeta)=\sqrt{2\pi}\langle\,\tilde u, (B_{z,\zeta}F^*\otimes I)
\tilde v_l\,\rangle 
$,
where $v_l$ is the sequence given by Lemma~\ref{lemvl}, and $B_{z,\zeta}$ is 
what one gets in \eqref{bzz} making $A=O(a)$. Since, for every 
$f\in{\mathcal S}(\re^2)$, 
\[
||\int[c_l(z,\zeta)-c(z,\zeta)]f(z,\zeta)\dz\dzeta||
\leq 
\]\[
||u||_{L^2}\cdot||B||\cdot||v-v_l||_{L^2}\cdot
\int|f(z,\zeta)|\,\dz\dzeta\,\longrightarrow\,0,
\]
as $l\to\infty$, it is enough to show that 
\[
\lim_{l\to\infty}\int[c_l(z,\zeta)-a(z,\zeta)]f(z,\zeta)\dz\dzeta\,=\,0.
\]

It follows from \eqref{cvm} that $B=O(b)$, for $
b(x,\xi)=(1+\partial_x)^2(1+\partial_\xi)^2a(x,\xi)$.
We then get $B_{z,\zeta}=O(b_{z,\zeta})$, for 
$b_{z,\zeta}(x,\xi)=b(x+z,\xi+\zeta)$. Hence, if $\varphi$\ and $\psi$\ 
belong to ${\mathcal S}(\re)$, then 
\[
[(B_{z,\zeta}F^{*}\otimes I)(\varphi\otimes\psi)](x,\eta)
\,=\,
\frac{1}{\sqrt{2\pi}}\int e^{i
x\xi}b(x+z,\xi+\zeta)\varphi(\xi)\psi(\eta)\dxi.
\]
Using that 
$v_l\in{\mathcal S}(\rn)\talg{\mathcal S}(\rn)$, we then get
\[
c_l(z,\zeta)=\int_{\re^{3}}
\overline{u(x,\xi)}e^{ix\xi}b(x+z,\xi+\zeta)v_l(\xi,\eta)
\dxi\dx\deta.
\]
By Lemma~\ref{cordes}, we then have
\[
\int[c_l(z,\zeta)-a(z,\zeta)]f(z,\zeta)\dz\dzeta=
\]
\[
\int_{\re^2}[\int_{\re^3}
\overline{u(x,\eta)}e^{ix\xi}b(x+z,\xi+\zeta)(v(\xi,\eta)-v_l(\xi,\eta))
\dxi\dx\deta]
f(z,\zeta)\dz\dzeta.
\]
Since $(x,\xi,\eta)\mapsto\overline{u(x,\eta)}(v(\xi,\eta)-v_l(\xi,\eta))$
belongs to $L^1(\re^3)$, we may interchange the order of integration and obtain that the above expression is bounded by
\[
\sup_{x,\xi}||b(x,\xi)||\cdot||f||_{L^1}\cdot
\int_{\re^{3}}|u(x,\eta)|\cdot|v(\xi,\eta)-v_l(\xi,\eta)|\,\dxi\dx\deta,
\]
which tends to zero, by Lemma~\ref{lemvl}, as we wanted. 

This proves that $S$ is a left-inverse for $O$ when $n=1$. We now comment on some of the changes needed to extend these definitions and
proof for arbitrary $n$. We have to replace $u$ and $v$, respectively, 
by $u_n(x,\eta)=u(x_1,\eta_1)\cdots u(x_n,\eta_n)$ and 
$v_n(\xi,\eta)=v(\xi_1,\eta_1)\cdots v(\xi_n,\eta_n)$. In
the definitions on $S$ and $c_l$, we replace $\sqrt{2\pi}$\ by 
$(2\pi)^{n/2}$, $\langle\cdot,\cdot\rangle$ denotes the inner product of
$E_{2n}$ and $F\in{\mathcal B}^*(E_n)$. The new $B_{z,\zeta}$ is defined by
\begin{equation}
\label{defbn}
B_{z,\zeta}
=[\prod_{j=1}^{n}(1+\partial_{z_{j}})^2(1+\partial_{\zeta_{j}})^2]A_{z,\zeta}.
\end{equation}
The integral in Lemma~\ref{cordes}\ is now an integral over
$\re^{3n}$\ and the equality in \eqref{hoc}\ holds for all
$(z,\zeta)\in\re^{2n}$. The integral in Lemma~\ref{lemvl} is also over
$\re^{3n}$, and $v_l$ belongs to ${\mathcal S}(\re^n)\talg{\mathcal S}(\re^n)$.

\section{Commutative Case}\label{s3}

In this section, we assume that $C$ is equal to $C(\Omega)$, the algebra of 
continuous functions on a Hausdorff compact topological space $\Omega$. For 
each $\lambda\in\Omega$ and each $f\in\sC$, we define 
$V_\lambda f\in{\mathcal S}(\rn)$ by
\[
(V_\lambda f)(x)\,=\,[f(x)](\lambda),\ \ x\in\rn.
\]
$V_\lambda$\ extends to a continuous linear mapping
$V_\lambda:E_n\longrightarrow L^2(\rn)$, with $||V_\lambda||\leq 1$.

\begin{lem}\label{lq}
Let there be given  $T\in\ben$, $f\in E_n$ and $\lambda\in\Omega$. If
$V_\lambda f=0$, then $V_\lambda Tf=0$.
\end{lem}

\pf
The equality $\langle V_\lambda g,V_\lambda g\rangle_{_{L^2(\rn)}}=\langle g, g\rangle
(\lambda)$ holds for all $g\in\sC$; hence also for all $g\in E_n$. We
then have:
\[
\langle V_\lambda Tf,V_\lambda Tf\rangle=\langle Tf,Tf\rangle(\lambda)=
\langle f,T^*Tf\rangle(\lambda)=|\langle f,T^*Tf\rangle(\lambda)|\leq\]\[
\sqrt{\langle f,f\rangle(\lambda)}\sqrt{\langle
T^*Tf,T^*Tf\rangle(\lambda)}=
\sqrt{\langle V_\lambda f,V_\lambda f\rangle_{_{L^2(\rn)}}}\sqrt{\langle V_\lambda T^*Tf,V_\lambda T^*Tf\rangle_{_{L^2(\rn)}}}.
\]
This implies our claim.\cqd

Given $\varphi\in{\mathcal S}(\rn)$, let $\tilde\varphi\in\sC$\ be defined by
$[\tilde\varphi(x)](\lambda)=\varphi(x)$, for all $\lambda\in\Omega$ and all $x\in\rn$. 
It is obvious that $V_\lambda\tilde\varphi=\varphi$. 
Given $T\in\ben$\ and $\lambda\in\Omega$, let $T_\lambda$ denote the unique linear
mapping defined by the requirement that the diagram
\[
\def\mapup#1{\Big\uparrow\rlap{$\vcenter{\hbox{$\scriptstyle#1$}}$}}
\def\mapdn#1{\Big\downarrow\rlap{$\vcenter{\hbox{$\scriptstyle#1$}}$}}
\begin{array}{ccc}
\sC&{\mathop{\longrightarrow}\limits^{T}}&E_n
\\&&
\\\mapdn{V_\lambda}&&\mapdn{V_\lambda}
\\&&
\\{\mathcal S}(\rn)&{\mathop{\longrightarrow}\limits^{T_\lambda}}&L^2(\rn)
\end{array}
\]
commutes. This is well defined by Lemma~\ref{lq} and because the left vertical
arrow in the above diagram is surjective. 

\begin{lem}\label{lr} For each $T\in\ben$ and each $\lambda\in\Omega$, $T_\lambda$ extends to a
bounded operator on $L^2(\rn)$. Moreover, we have
\begin{equation}\label{desr}
||T||\ =\ \sup\,\{||T_\lambda||;\,\lambda\in\Omega\}.
\end{equation}
\end{lem}
\pf
Given $\varphi\in{\mathcal S}(\rn)$, let $\tilde\varphi$ denote the element of
$E_n$ defined after Lemma~\ref{lq}. We have:
\[
||T_\lambda\varphi||_{_{L^2(\rn)}}=||V_\lambda T\tilde\varphi||_{_{L^2(\rn)}}\leq ||T\tilde\varphi||_2\leq||T||\cdot||\tilde\varphi||_2=
||T||\cdot||\varphi||_{_{L^2(\rn)}}.
\]
This implies that $T_\lambda$ extends to a bounded operator on $L^2(\rn)$ with
norm bounded by $||T||$. 

Let $M$ denote the right-hand side of (\ref{desr}). For each
$\lambda\in\Omega$ and each $f\in\sC$, using Lemma~\ref{lq} and the first
statement in its proof, we get:
\[
|\langle Tf,Tf\rangle(\lambda)|=|\langle V_\lambda Tf,V_\lambda Tf\rangle_{_{L^2(\rn)}}|
=\]\[|\langle T_\lambda V_\lambda f,T_\lambda V_\lambda
f\rangle_{_{L^2(\rn)}}|\leq
||T_\lambda||\cdot||V_\lambda f||_{_{L^2(\rn)}}\leq M||f||_2.
\]
Taking the supremum in $\lambda$ on the left, we get $||Tf||_2\leq M||f||_2$. \cqd

Our goal in this Section is to prove that the mapping $S$ defined in the
previous section is injective for $C=C(\Omega)$. This will finish the proof of
Theorem~\ref{t}. 

Given $A\in\h$ such that $SA=0$, we want to show that $A=0$. In view of
the following lemma, it suffices to show that $B=0$, where $B=B_{0,0}$
($B_{z,\zeta}$ as defined on \eqref{defbn}). Lemma~\ref{ls}\ is \cite[Proposition 3.1]{book} when $C=\C$. The same proof works
for any C$^*$-algebra $C$.

\begin{lem}\label{ls}
If $Y\in\h$, $Y_{z,\zeta}=T_{-z}M_{-\zeta}YM_{\zeta}T_z$ $(z,\zeta\in\rn)$, and either $(1+\partial_{z_{j}})Y_{z,\zeta}\equiv 0$ or
$(1+\partial_{\zeta_{j}})Y_{z,\zeta}\equiv 0$ for some $j$, then $Y=0$
\end{lem}

By Lemma~\ref{lr}, in order to prove that $B=0$, it suffices to show that
$B_\lambda=0$ for each $\lambda\in\Omega$. 
For $z$ and $\zeta$ in $\re^n$, define $E_{z,\zeta}=M_\zeta T_z$. We then have
$B_{z,\zeta}=E_{z,\zeta}^*BE_{z,\zeta}$. Using that
$E_{z,\zeta}F^*=e^{iz\cdot\zeta}F^*E_{\zeta,-z}$, we may rewrite equation
$SA=0$ as 
\[e^{iz\cdot\zeta}
\langle(E_{z,\zeta}\otimes I)\tilde{u}_n,(BF^*E_{\zeta,-z}\otimes
I)\tilde{v}_n\rangle=0,\ \ \mbox{for all}\ \ (z,\zeta).
\]
Evaluating  this equation at $\lambda$ gives:
\begin{equation}
\label{agora}
e^{iz\zeta}\langle(E_{z,\zeta}\otimes I)u_n,(B_\lambda F^*E_{\zeta,-z}\otimes
I)v_n\rangle_{_{L^2(\re^n)}}=0,\ \ \mbox{for all}\ \ (z,\zeta).
\end{equation}

For a fixed $\varphi\in C_c^\infty(\re^{2n})$ to be chosen soon, and for each bounded operator
$D$ on $L^2(\re^n)$, define 
\begin{equation}
\label{defXi}
\Xi(D)\,=\,\int\varphi(z,\zeta)e^{iz\zeta}\langle(E_{z,\zeta}\otimes I)u_n,(DF^*E_{\zeta,-z}\otimes
I)v_n\rangle_{_{L^2(\re^n)}}\dz\dzeta.
\end{equation}
In case $D$ is finite-rank, and hence we may take $b^1,\cdots,b^k,c^1,\cdots,c^k$ in $L^2(\re^n)$ such that, for all $f\in
L^2(\re^n)$, 
\[
D F^*f=\sum_{j=1}^kb^j\langle c^j,f\rangle_{_{L^2(\re^n)}},
\]
we have: $\Xi(D)=$
\[
\sum_{j=1}^k\int\!\!\!\int b^j(x)\bar{c}^j(\xi)\int\!\!\!\int\!\!\!\int e^{iz\cdot\zeta}
\varphi(z,\zeta)e^{-ix\cdot\zeta}\bar{u}_n(x-z,\eta)e^{-iz\cdot\xi}v_n(\xi-\zeta,\eta)\dz\dzeta\deta\dxi\dx.
\]
Making the change of variables $x-z=z^\prime$, $\xi-\zeta=\zeta^\prime$\ on
the inner triple integral above, we get:
\begin{equation}
\label{quasi}
\Xi(D)=\sum_{j=1}^k\int\!\!\!\int b^j(x)\bar{c}^j(\xi)e^{-ix\cdot\xi}\int\!\!\!\int\!\!\!\int e^{iz\cdot\zeta}
\varphi(x-z,\xi-\zeta)\bar{u}_n(z,\eta)v_n(\zeta,\eta)\dz\dzeta\deta\dxi\dx.
\end{equation}
For arbitrary $\chi$ and $\psi$ in $C_c^\infty(\re^n)$, let $\varphi$ be defined by 
\[
(1+\partial_x)^2(1+\partial_\xi)^2[e^{ix\xi}\bar\chi(-x)\psi(-\xi)]=\varphi^\sharp(x,\xi),\
\ \ \varphi(x,\xi)=\varphi^\sharp(-x,-\xi). 
\]
Using the higher dimensional version of Lemma~\ref{cordes} mentioned at the
end of Section~\ref{s2}, the right side of \eqref{quasi} becomes:
\[
\sum_{j=1}^k\int\!\!\!\int
b^j(x)\bar{c}^j(\xi)\bar\chi(x)\psi(\xi)\dx\dxi=\langle\chi,D
F^*\psi\rangle_{_{L^2(\re^n)}}.
\]
This shows that, for this choice of $\varphi$, 
\begin{equation}
\label{posto}
\Xi(D)=\langle\chi,DF^*\psi\rangle_{_{L^2(\re^n)}},
\end{equation} 
whenever $D$ has finite rank. 

Let $\{\phi_1,\phi_2,\cdots\}$ be an orthonormal basis of ${L^2(\re^n)}$. For
each positive integer $j$, let $P_j$ denote the orthogonal projection onto the spam of 
$\{\phi_1,\cdots,\phi_j\}$. Cordes proved (\cite[Chapter 8]{book}, between equations (3.27) and
(3.29)) that, for any bounded operator $T$ on $L^2(\re^n)$, one has
$\lim_{j\to\infty}\Xi(P_jTP_j)\ =\ \Xi(T)$. Applying this to $T=B_\lambda$ and
using \eqref{posto}, we get
\[ 
\Xi(B_\lambda)=\lim_{j\to\infty}\Xi(P_jB_\lambda P_j)=
\lim_{j\to\infty}\langle\chi,P_jB_\lambda P_jF^*\psi\rangle_{_{L^2(\re^n)}}
\ =\ \langle\chi,B_\lambda F^*\psi\rangle_{_{L^2(\re^n)}}.
\]
By \eqref{agora}, the left-hand side of this equality vanishes. Since $\chi$
and $\psi$ are arbitrary test functions, this shows that $B_\lambda=0$. This finishes the proof of Theorem~\ref{t} (recall our remarks
before and after the statement of Lemma~\ref{ls}).

\section{Rieffel's Conjecture}\label{s4}

Given a skew-symmetric $n\times n$ matrix $J$ and $F\in\bcn$ (i.e., $F:\rn\to C$ is smooth and, together with all its 
derivatives, is bounded), let us denote by $L_F$ the \psd\ \op\
$a(x,D)\in\ben$ with symbol $a(x,\xi)=F(x+J\xi)$. At the end of Chapter~4 in
\cite{R}, Rieffel made a conjecture that may be rephrased as follows: any
operator $A\in\ben$ that is Heisenberg-smooth and commutes with every operator
of the form $R_G=b(x,D)$, where $b(x,\xi)=G(x-J\xi)$ with $G\in\bcn$, is of the
form $A=L_F$ for some $F\in\bcn$. 

Using Cordes characterization of the Heisenberg-smooth operators in the scalar
case, we have shown \cite{MM} that Rieffel's conjecture is true when
$C=\C$. It has been further proven by the second-named author \cite{M} that Rieffel's
conjecture is true for any separable $C^*$-algebra $C$ for which the operator $O$
defined in \eqref{o} is a bijection. Under this assumption, a result actually stronger than what was
conjectured by Rieffel was proven in \cite[Theorem 3.5]{M}: To get $A=L_F$ for
some $F\in\bcn$, one only needs to require that a given $A\in\ben$ is
``translation-smooth'' (i.e., the mapping $\rn\ni z\mapsto  T_{-z}AT_z\in\ben$
is smooth) and commutes with every $R_G$ with $G\in\sC$. Combining this result with our Theorem~\ref{t}, we then get:

\begin{thm} Let $C$ be a unital commutative separable C$^*$-algebra. If a given
$A\in\ben$ is translation-smooth and commutes with every $R_G$, $G\in\sC$,
then $A=L_F$ for some $F\in\bcn$.
\end{thm}

\section*{Acknowledgements} Severino Melo was partially supported by the
Brazilian agency CNPq (Processo 306214/2003-2), and Marcela Merklen had a posdoc
position sponsored by CAPES-PRODOC. We thank Ricardo
Bianconi, Ruy Exel and Jorge Hounie for several helpful conversations.

\vskip1.0cm

{\footnotesize 

\noindent
Instituto de  Matem\'atica e Estat\'{\i}stica\\
Universidade de S\~ao Paulo\\
Caixa Postal 66281 \\
05311-970 S\~ao Paulo, Brazil. 

\noindent
Email: toscano@ime.usp.br, marcela@ime.usp.br}


\begin{thebibliography}{99}

\bibitem{CV} {\sc A. P. Calder\'on \&\ R. Vaillancourt}, {\em On the 
boundedness of pseudo-differential operators}; J. Math. Soc. Japan {\bf 23} (1971), 
374-378.

\bibitem{C} {\sc H. O. Cordes}, On \psd\ \ops\ and smoothness of special 
Lie-group representations, Manuscripta Math. {\bf 28} (1979), 51-69.

\bibitem{book} {\sc H. O. Cordes}, The technique of Pseudodifferential 
Operators, London Mathematical Society Lecture Note Series {\bf 202},
Cambridge Univesity Press, Cambridge, 1995.

\bibitem{H} {\sc I. L. Hwang}, {\em The $L^2$-boundedness of
pseudodifferential operators}, Trans. Amer. Math. Soc. {\bf 302}-1 (1987),
55-76. 

\bibitem{L} {\sc C. Lance}, Hilbert C$^*$-modules - A toolkit for
operator algebraists, London Mathematical Society Lecture Note Series 
{\bf 210}, Cambridge University Press, Cambridge, 1995. 

\bibitem{MM} {\sc S. T. Melo \&\ M. I. Merklen}, {\em On a conjectured
Beals-Cordes-type characterization}, Proc. Amer. Math. Soc. {\bf 130}-7
(2002), 1997-2000.

\bibitem{M} {\sc M. I. Merklen}, {\em Boundedness of Pseudodifferential 
Operators of C$^*$-Algebra-Valued Symbol}; Proc. Roy. Soc. Edinburgh Sect. A
{\bf 135}-6 (2005), 1279-1286.

\bibitem{Mt} {\sc M. I. Merklen}, Resultados motivados por uma 
caracteriza\c c\~ao de operadores pseudo-diferenciais conjecturada por
Rieffel, Tese de Doutorado, Universidade de S\~ao Paulo, 2002, 
http://arxiv.org/abs/math.OA/0309464. 

\bibitem{Ra} {\sc M. Rieffel}, {\em Induced Representations of C$^*$-Algebras}, 
Advances in Math. {\bf 13} (1974), 176-257.

\bibitem{R} {\sc M. Rieffel}, Deformation Quantization for Actions of 
${\mathbb R}^d,$\ Memoirs of the American Mathematical Society {\bf 506}, 1993.

\bibitem{S} {\sc J. Seiler}, {\em Continuity of edge and corner
pseudodifferential operators}, Math. Nachr. {\bf 205} (1999), 163-182.

\end{thebibliography}
\end{document}